\documentclass[]{interact}

\usepackage{epstopdf}
\usepackage{subfigure}

\theoremstyle{plain}
\newtheorem{theorem}{Theorem}[section]
\newtheorem{lemma}[theorem]{Lemma}

\usepackage{mathtools}
\DeclarePairedDelimiter{\ceil}{\lceil}{\rceil}
\theoremstyle{definition}

\theoremstyle{remark}

\usepackage{algorithm}
\usepackage[noend]{algpseudocode}
\renewcommand{\algorithmiccomment}[1]{\bgroup\hfill//~#1\egroup}
\usepackage{amsmath}
\usepackage{amsthm}
\usepackage{amssymb}
\usepackage{graphicx}
\usepackage{adjustbox}
\usepackage{changepage}
\usepackage{setspace}
\usepackage[table]{xcolor}
\usepackage{diagbox}
\usepackage{xcolor}
\usepackage{tikz,pdfcomment}
\usepackage{pgfplots}
\usepackage{multicol}
\usepackage{multirow}
\usepackage{pdflscape}
\usepackage{hyperref}
\usepackage{comment}
\usepackage{verbatim}

\usepackage{pgfplots, pgfplotstable}

\definecolor{A}{HTML}{D7191C}
\definecolor{B}{HTML}{D7198C}
\definecolor{C}{HTML}{FDAE61}
\definecolor{D}{HTML}{008000}
\definecolor{E}{HTML}{FE6F5E}
\definecolor{F}{HTML}{ABDDA4}
\definecolor{G}{HTML}{BEBEBE}
\definecolor{H}{HTML}{808080}
\let\svthefootnote\thefootnote
\newcommand\freefootnote[1]{%
  \let\thefootnote\relax%
  \footnotetext{#1}%
  \let\thefootnote\svthefootnote%
}

\begin{document}
\title{Weighted tardiness minimization for unrelated machines with sequence-dependent and resource-constrained setups}

\author{
\name{Ioannis Avgerinos \textsuperscript{a}, Ioannis Mourtos\textsuperscript{a}, Stavros Vatikiotis\textsuperscript{a} and Georgios Zois\textsuperscript{a}\thanks{CONTACT Vatikiotis Stavros. Email: stvatikiotis@aueb.gr. Department of Management Science and Technology, Athens University of Economics and Business, Athens 104 34, Greece}} 
\affil{\textsuperscript{a} ELTRUN Research Lab, Department of Management Science and Technology, Athens University of Economics and Business, Athens 104 34, Greece}
}

\maketitle

\begin{abstract}
\indent Motivated by the need of quick job (re-)scheduling, we examine an elaborate scheduling environment under the objective of total weighted tardiness minimization. The examined problem variant moves well beyond existing literature, as it considers unrelated machines, sequence-dependent and machine-dependent setup times and a renewable resource constraint on the number of simultaneous setups. For this variant, we provide a relaxed MILP to calculate lower bounds, thus estimating a worst-case optimality gap. As a fast exact approach appears not plausible for instances of practical importance, we extend known (meta-)heuristics to deal with the problem at hand, coupling them with a Constraint Programming (CP) component - vital to guarantee the non-violation of the problem's constraints – which optimally allocates resources with respect to tardiness minimization. The validity and versatility of employing different (meta-)heuristics exploiting a relaxed MILP as a quality measure is revealed by our extensive experimental study, which shows that the methods deployed have complementary strengths depending on the instance parameters. Since the problem description has been obtained from a textile manufacturer where jobs of diverse size arrive continuously under tight deadlines, we also discuss the practical impact of our approach in terms of both tardiness decrease and broader managerial insights.
\end{abstract}



\section{The problem} \label{Section:Problem}
 
\textbf{Background.} The common occurrence of scheduling problems in both logistics and manufacturing has motivated the integration of optimisation methods into the decision-making process. Although most scheduling problems and their computational complexity have been defined several decades ago (e.g., \cite{Len77}), the transition to Industry 4.0 has increased, simultaneously, the availability of credible production data (e.g., processing times or setup times), the diversity of features appearing in practical cases, and the strength of computational tools. It is therefore not surprising that the interest in production scheduling deepens, as indicatively discussed in two recent reviews \cite{Jiang22, Par20}.

The numerous variants of scheduling in the literature match the variety of scheduling problems formulated in practice. These variants are distinguished by the number of the available machines, the differences in their speeds, the sequence-dependency of setup times, the eligibility of machines per job, the resource consumption rate, the renewability of the resources, and, of course, the objective. Although the most common \cite{Lawler93} objective remains the minimisation of the \emph{makespan}, i.e., the maximum completion time over all jobs in a feasible schedule, modern manufacturing is continuously facing disruptions and tight deadlines, thus imposing a research shift towards due-time related optimisation criteria. Perhaps the most realistic such criterion is the delay beyond a deadline (i.e., the \emph{tardiness}) weighted by a job-specific value to differentiate either the job importance or actual financial penalties associated with breaching a deadline \cite{Len77}. 

Indeed, the optimisation community is recently focusing on tardiness-related objectives \cite{Kuh16}. Still, the cases examined up to now omit critical shop-floor realities such as sequence-dependence, unrelated parallel machine environments and limited resource availability. Also, optimising tardiness instead of makespan appears to have a negative impact on the performance of exact methods, hence the challenge of solving instances of realistic size even to (provably) half-optimality. This is where our contribution comes forward: we modify state-of-the-art (meta-)heuristics to minimise total weighted tardiness in an unrelated parallel machine environment with resources and sequence-dependencies. To establish valid optimality gaps, we propose a MILP that obtains lower (i.e., optimistic or so-called `dual') bounds. We show that there are no globally preferable methods, as the performance of the presented (meta-)heuristics is highly dependent on the structure or the scale of the instance examined. Our combined methods find competitive solutions, i.e., at least half-optimal for instances where our dual-bounding MILP is solved to optimality within one hour (for up to $50$ jobs and $10$ machines). The competitive behaviour is sustained for benchmark instances of up to $200$ jobs and $20$ machines and also in large-sized real instances.

\textbf{Literature Review.} The simplest variant of a parallel machine scheduling concerns identical machines, meaning that each machine requires an identical amount of time to process the same job. The reduced complexity of this variant allows the identification of properties of optimal schedules, as shown by \cite{Yal00} for the minimisation of total weighted completion times, and by \cite{Salah} for the minimisation of total weighted tardiness. Since the industrial settings usually involve non-identical machines, there are several recent implementations of meta-heuristic and exact methods on such settings, e.g., \cite{Kim20a, Lee18, Pinedo97,Oz18}, to name a few. However, none of these methods examines tardiness in such an elaborate setting as ours.

Two interesting variants are uniform and unrelated machines. Regarding the former, an elementary study on parallel uniform machines, with a due-time related objective, appears in \cite{Emmons87}. The extended variant of the minimisation of makespan on uniform machines problem of \cite{Lee21} includes additional features met in practice, such as sequence-dependencies or job splitting. A Variable Neighbourhood Search approach is implemented on unrelated machines to minimise total weighted tardiness by \cite{Bilyk}. An exact mathematical model and a heuristic approach for the same objective are presented by \cite{Bek19}. The work of \cite{Chen09, Lee13} regarding the minimisation of total single-weighted tardiness is also motivating. Recently, a set of exact and (meta)heuristic methods, including local search operators and simulated annealing (SA) algorithms, for a bi-objective scheduling problem, is described by \cite{Mos21}. Two recently implemented methods rely on Variable Neighbourhood Descent \cite{Liao22} and GRASP algorithms \cite{Iori22}. Although all the above examine production settings without any resource constraints, they motivate the use of (meta-)heuristics on unrelated machine scheduling with sequence-dependency to minimise (total weighted) tardiness.

Regarding exact methods on due-time related objectives, we observe a major gap in relevant literature. Existing MILP formulations solve small-sized instances, although, for example, the real-life instances we tackle assume $50$ jobs on average, processed by $5-10$ machines. This is of no surprise, since, a typical MILP introduces `big-M' constraints to ensure the continuity of sequences, thus yielding very weak lower bounds. In addition, the tardiness $T_{j}$ of a job $j$ is formulated as the nonlinear equation  $T_{j} = max\{0, C_{j} - d_{j}\}$, where $C_{j}$ is the completion time and $d_{j}$ the deadline of $j$; the (trivial) linearisation of this expression deteriorates the lower bound significantly. The exact approach of \cite{Pes22} employs column-generation approach for a time-indexed MILP formulation tackling tardiness on a single-machine but for batch scheduling. The progress on tardiness is in strike contrast to the makespan, which is better handled for elaborate scheduling problems through exact methods \cite{IJPR} or MILP-based heuristics \cite{Kim20b}. Here we make use of exact methods in two manners: first, we use a relaxed MILP formulation that provides good lower bounds but of course requires considerable computational time; secondly, we handle resource constraints in an exact manner by finding the optimal resource allocation given a sequencing of jobs to machines.

Indeed, the most challenging features of parallel machine scheduling arise from the restrictions imposed by available resources, which are consumed or occupied when a scheduled task is processed. Constraint Programming (CP) appears to perform better under such restrictions, through \emph{global constraints} that handle resource-allocation more easily than, e.g., time-indexed MILP models. \cite{Hook07} implemented a generic Logic-Based Benders Decomposition (LBBD) for scheduling under various objectives, by formulating a scheduling MILP master problem, assuming infinite resources, then obtaining primal solutions by a CP subproblem. A similar decomposition-based exact method, regarding the minimisation of makespan for a setting of unrelated machines, is presented by \cite{Peyro20}. To the contrary, \cite{Yun22} opted for a monolithic CP formulation to solve their multi-resource-constrained variant of unrelated machines scheduling to minimise makespan. Meta-heuristic methods are also of value, as \cite{Pin20} being an indicative example of integrating resource constraints into simulated annealing (SA) and local search algorithms. Our LBBD scheme \cite{MIM} minimises total (unweighted) tardiness for unrelated machines but for instances where setups are of small duration and sequence-dependency is not as intensive as in the instances examined here.

To sum up, Table \ref{tab:review} lists the scheduling problems presented recently, distinguished by machine environment, sequence-dependency, resources availability, objectives (e.g., makespan, tardiness) and solution methodologies (e.g., exact methods, heuristics, meta-heuristics). The families of scheduling problems are denoted by the standard three-field notation  $\alpha|\beta|\gamma$ of \cite{Grah79}, where $\alpha$ is refers to the machine environment, i.e., 1 for single machine, $P$ for identical machines, $Q$ for uniform machines, and $R$ for unrelated machines, $\beta$ refers to the job properties e.g. deadlines, sequence-dependencies ($ s_{jk}$), precedence constraints ($prec$), release dates ($r_j$),  and $\gamma$ indicates the objective e.g., $C_{max}$ for makespan, $\sum_{j}C_{j}$ for total completion times, $\sum_{j}w_{j}T_{j}$ for total weighted tardiness.  

\begin{table}[H]
\centering
\caption{Review of Scheduling Problems (2017-2022)}
\resizebox{\textwidth}{!}{
\begin{tabular}{ccccccc}
\hline
Reference        & Objective                  & Machines & Sequence-Dependency & Resources  & Exact Method & (Meta)heuristics \\ \hline
Eroglu \& Ozmutly, 2017\cite{Eroglou17} & $C_{max}$                  & R        & \checkmark          & -          & -            & GA               \\
Hamzadayi \& Yildiz, 2017 \cite{Hamz17}    & $C_{max}$                  & P        & \checkmark          & \checkmark & MILP         & GA, SA           \\
Lee, 2018 \cite{Lee18}     & $\sum_{j}T_{j}$            & P        & \checkmark          & -          & MILP         & ATCS             \\
Ozturk \& Chu, 2018 \cite{Oz18}      & $\sum_{j}T_{j}$            & P        & -                   & -          & MILP         & GA               \\
Bektur \& Saraç, 2019 \cite{Bek19}     & $\sum_{j}w_{j} T_{j}$ & R        & \checkmark          & - & MILP         & ATCS, SA, TS     \\
Kim et al., 2020 \cite{Kim20a}    & $\sum_{j}T_{j}$            & P        & \checkmark          & -          & MILP         & GA, SA           \\
Kim \& Kim, 2020 \cite{Kim20b}    & $C_{max}$                  & P        & \checkmark          & -          & MILP         & GA               \\
Peyro, 2020 \cite{Peyro20}   & $C_{max}$                  & R        & \checkmark          & \checkmark & D            & -                \\
Pinheiro et al., 2020 \cite{Pin20}     & $\sum_{j}T_{j}$            & R        & \checkmark          & \checkmark & MILP         & SA, LS           \\
Kim \& Lee, 2021 \cite{Lee21}     & $C_{max}$                  & Q        & \checkmark          & \checkmark & MILP         & Other            \\
Moser et al., 2021 \cite{Mos21}     & $C_{max}$, $\sum_{j}T_{j}$ & R        & \checkmark          & -          & MILP         & LS, SA, Other    \\
Avgerinos et al., 2022 \cite{IJPR}      & $C_{max}$                  & R        & \checkmark          & \checkmark & LBBD         & Other            \\
Avgerinos et al., 2022 \cite{MIM}       & $\sum_{j}T_{j}$            & R        & \checkmark          & \checkmark & LBBD         & -                \\
Iori et al., 2022 \cite{Iori22}    & $\sum_{j}w_{j} T_{j}$ & R        & -                   & -          & -            & GRASP            \\
Liao et al., 2022 \cite{Liao22}    & $C_{max}$                  & R        & \checkmark          & -          & MILP         & VNS              \\
Pessoa et al., 2022 \cite{Pes22}     & $\sum_{j}w_{j} T_{j}$ & 1        & -                   & -          & CG           & -                \\
Yunusoglu \& Yildiz, 2022 \cite{Yun22}     & $C_{max}$                  & R        & \checkmark          & \checkmark & CP           & -                  \\
Current Work    & $\sum_{j}w_{j} T_{j}$    & R     & \checkmark            & \checkmark & MILP & GA, SA, ATCS              \\ \hline
\end{tabular}
}
{\tiny GA: Genetic Algorithm, MILP: Mixed-Integer Linear Programming, SA: Simulated Annealing, ATCS: Apparent Tardiness Cost with Setups, TS: Tabu Search, D: Decomposition, LS: Local Search, GRASP: Greedy Randomized Adaptive Search Procedure, VNS: Variable Neighbourhood Search, CG: Column Generation, CP: Constraint Programming}
\label{tab:review}
\end{table}

\noindent\textbf{Our contribution.} We advance the current literature by optimising total weighted tardiness under unrelated parallel machines, sequence and machine-dependent setup times and a renewable constraint on the total number of simultaneous setups. Indicatively, we extend the work of \cite{Bek19} by assuming a resource in multiple copies (instead of a single common resource or ``server") and by omitting restrictions like machine eligibility. In that sense, we also deal with the entire family of problems in \cite{Peyro20} but for the objective of total weighted tardiness. 

To evaluate the (worst-case) quality of a solution, we provide a MILP, formulating a relaxed version of the problem. To our knowledge, this the first lower-bounding MILP in the presence of unrelated machines with non-zero setups and for weighted tardiness minimisation (see \cite{Liaw03} for a bounding argument that is the least distant to ours). 

To find solutions of good quality, in the absence of scalable exact methods, we modify state-of-the-art (meta-)heuristics to deal with tardiness and the existence of resources to be allocated. Our core modification is a CP model which, given a sequencing of jobs per machine, allocates resources optimally with respect to weighted tardiness minimisation. 

In addition, we experiment on novel benchmark instances, as our problem variant has not been studied elsewhere, and find good-quality solutions for up to $200$ jobs and $20$ machines within an hour. For instances up to $50$ jobs and $10$ machines, we  allow for an additional 1-hour to compute an estimate on the optimality gap through our MILP. Our experimentation shows a dependency between the size of the instance and the efficiency of each algorithm. Specifically, the proposed Genetic Algorithm (GA) outperforms the ATCS algorithm for instance sizes of up to 100 jobs, while the SA algorithm, even given a warm-start solution from the aforementioned methods, does not justify the computational time consumed. However, when the number of jobs exceeds $100$, the GA is not able to compete with the ATCS: the explanation, probably, lies on the fact that the solution space grows exponentially, thus the GA has not enough time to explore it in reasonable time. 

To investigate our contribution in practical terms, the proposed algorithms are employed on real datasets from textile manufacturing. Here, apart from creating schedules of good quality in instances of varying size, our work serves as a simulation tool to determine potential improvements on the production process and evaluate possible investments on behalf of the plant managers.

\noindent\textbf{Outline.} This remainder of this paper is structured as follows. In Section 2, after a formal description of the problem, we present a novel MILP formulation that computes valid lower bounds. This relaxation is an exact solving method for the sequence-independent variant of the problem, as our approach replaces sequence-dependent times with `optimistic' estimations. In Section 3, we show a CP module for the resource allocation part of the problem and, then, modify three (meta-)heuristics to make them applicable on our variant. In Section 4, we implement and evaluate our algorithms on two sets of instances. The first one is constructed by a generator of random instances, considering multiple variations for all involved parameters and relying upon literature generators to the largest possible extent. The second one is obtained by a real industry operator, originating from the textile industry. We conclude in Section 5 by discussing managerial insights and listing  ideas for further research.  
\section{Dual bounds through MILP}

Let $J$ be a set of jobs. Each job $j\in J$ is featured with a weight value $w_{j}$ and a deadline $d_{j}$, meaning that each time unit of delay implies a penalty cost $w_{j}$. A set of machines $M$ is assigned to perform jobs $J$. Each machine $m\in M$ can process one job at a time and multiple machines can process different jobs simultaneously. Machines $M$ are unrelated, therefore a processing time $p_{jm}$ is attributed to each job $j\in J$. Each job, before being processed, should be set up to the assigned machine. The respective setup times $s_{ijm}$ depend on both the machine and the sequence, meaning that they are determined by the preceding job $i\in J\setminus\{j\}$ in the same machine $m\in M$. Each setup task occupies a unit of working resources $WR$, therefore no more than $WR$ setup tasks can be conducted simultaneously. Resources are renewable, meaning that once a unit is available again, it can be occupied by a new setup task. We aim at computing a schedule at the minimum total weighted tardiness. Table \ref{Tbl2} lists all relevant notation.

As mentioned above, to evaluate the divergence of the solutions of the proposed heuristic methods from optimality (i.e., the primal bounds), we opt for the construction of a valid MILP relaxation. This is obtained by replacing each sequence-dependent setup time with the minimum setup time for the respective job across all machines and precedent jobs. In that manner we omit sequence-dependent constraints, which are mainly responsible for the weak lower bounds of MILP formulations, and therefore achieve faster convergence to the optimal solution by an MILP solver.

\begin{table}[H]
    \fontsize{6}{12}\selectfont
    \tbl{Model Parameters}
    {\resizebox{\textwidth}{!}{
    \begin{tabular}[t]{ll}
    \hline
    \textbf{Sets} &\\
    \hline
    $J$         & Jobs \\
    $M$         & Machines \\
    \hline
    \textbf{Parameters} &\\
    \hline
    $p_{jm}$    & Processing time of job $j\in J$ to machine $m\in M$ \\
    $s_{ijm}$   & Setup time of job $j\in J$ succeeding job $i\in J\setminus \{j\}$ to $m\in M$ \\
    $s^{0}$     & Setup time of the first job of any machine \\ 
    $d_{j}$     & Deadline of job $j\in J$ \\ 
    $w_{j}$     & Weight of job $j\in J$ \\ 
    $WR$        & Number of re-usable working resources \\
    $s^{\leq}_{jm}$ &A lower bound of setup time of job $j\in J$ to machine $M\in M$; $s^{\leq}_{jm} = min_{i\in J\setminus \{j\}}\{s_{ijm}\}$ \\
    $t_{max}$   & The largest distinct time instance of the schedule \\
    $V$         & A large real number\\
    \hline
    \textbf{Variables} &\\
    \hline
    $x_{ijm}\in \{0,1\}$    & 1 if job $j\in J$ is assigned to slot $i\in J$ of machine $m\in M$, 0 otherwise \\
    $y_{jm}\in \{0, 1\}$    & 1 if job $j\in J$ is assigned to machine $m\in M$, 0 otherwise \\
    $f_{im}\in \{0, 1\}$    & 1 if slot $i\in J$ is the first occupied to machine $m\in M$, 0 otherwise \\
    $W_{ijmt}\in \{0, 1\}$  & 1 if job $j\in J$, assigned to slot $i\in J$ of machine $m\in M$, is delayed by $t\in \{0, ..., t_{max}\}$ distinct time instance, 1 otherwise \\
    $C_{im}\geq 0$          & Completion time of job or slot $i\in J$ to machine $m\in M$\\
    $P_{im}\geq 0$          & Processing time of slot $i\in J$ \\
    $S_{im}\geq 0$          & Setup time of slot $i\in J$ \\
    $D_{im}\geq 0$          & Deadline of slot $i\in J$ \\
    $T_{im}\geq 0$          & Tardiness of slot $i\in J$ \\
    $\lambda_{m_{i}}:$ \texttt{Interval}        & Interval variable, indicating the start and end time of the setup process of the $i^{th}$ job assigned to machine $m\in M$ \\ \hline
    \end{tabular}}}
    \label{Tbl2}
\end{table}
 
Let $\mathcal{P}$ be the original optimisation problem that includes setup times $s_{ijm}$ that depend on both the machine and the sequence. We construct a relaxation $\mathcal{R}$ of $\mathcal{P}$ as follows. For each job $j\in J$ and machine $m\in M$, let the original setup times values $s_{ijm}$ be replaced with $s^{\leq}_{jm} = min_{i\in J\setminus\{j\}}\{s_{ijm}\}$ (i.e., the minimum setup time for all candidate preceding jobs $i\in J\setminus\{j\}$). Let also $WR = |M|$, meaning in practice that the number of resources is unlimited. We observe that, for this case, the sequence-dependencies have been omitted, as the setup times are identical for any preceding job in the same machine. 

This relaxation, concerning unweighted tardiness, shows up as the master problem of the LBBD in \cite{MIM}. Our computational work there showed that, for instances which in which sequence-dependent tasks do not vary a lot, a reasonable amount of time is adequate to reach near-optimal solutions, but for total (i.e., unweighted) tardiness. Extending this approach to weighted tardiness  cannot solve instances of more than $10$ jobs to optimality, mainly because it struggles to improve lower bounds. However, we can still use it to obtain valid dual (aka lower) bounds.

To exploit the exclusion of sequence-dependencies from the presented relaxation $\mathcal{R}$, we opt for a `job-to-slot' assignment. That is, we assume that each machine has $|J|$ slots, each occupied by one job at most, hence each job is assigned to exactly one slot of one machine. 

We use the following binary variables: $y_{jm}$ indicates whether job $j\in J$ is assigned to machine $m\in M$ and $x_{ijm}$ is equal to 1 if job $j\in J$ is assigned to slot $i\in J$ of machine $m\in M$. Variables $f_{im}$ indicate whether $i\in J$ is the first occupied slot of machine $m\in M$. Continuous variables define the values of processing time, setup time, deadline, completion time and tardiness of slot $i\in J$ of machine $m\in M$ ($P_{im}$, $S_{im}$, $D_{im}$, $C_{im}$, $T_{im}$ respectively). Last, binary variables $W_{ijmt}$ are equal to 1 if job $j$, assigned to slot $i$ of machine $m$, is delayed by $t$ discrete time instances, and 0 otherwise. 

The relaxation $\mathcal{R}$ is as follows.

    \begingroup
    \scriptsize \selectfont
    \begin{flalign}
        \mathcal{R}: && \notag &&\\
        \text{min } & z\geq \sum_{i\in J}\sum_{j\in J}\sum_{m\in M}\sum_{t=0}^{t_{max}}w_{j}\cdot t\cdot W_{ijmt} \label{eq:M_obj} &&\\
        &\sum_{m\in M}y_{jm} = 1 && \forall j\in J \label{eq:m1} &&\\
        &\sum_{i\in J}x_{ijm} = y_{jm} && \forall j\in J, m\in M \label{eq:m2} &&\\
        &\sum_{j\in J}x_{ijm} \leq 1 && \forall i\in J, m\in M \label{eq:m3} &&\\
        &\sum_{i\in J}f_{im} = 1 && \forall m\in M \label{eq:m4} &&\\
        &\sum_{j\in J}x_{ijm} \geq \sum_{j\in J}x_{i-1jm} && \forall i\in J\setminus\{0\}, m\in M  \label{eq:m5} &&\\
        &\sum_{j\in J}x_{i-1jm} + f_{im} \leq 1 && \forall i\in J\setminus\{0\}, m\in M \label{eq:m6} &&\\
        &C_{im} = C_{i-1m} + P_{im} + S_{im} && \forall i\in J\setminus\{0\}, m\in M \label{eq:m7} &&\\
        &D_{im} = \sum_{j\in J}d_{j}\cdot x_{ijm} && \forall i\in J, m\in M \label{eq:m8} &&\\
        &T_{im} \geq C_{im} - D_{im} && \forall i\in J, m\in M \label{eq:m9} &&\\
        &P_{im} = \sum_{j\in J}p_{jm}\cdot x_{ijm} && \forall i\in J, m\in M \label{eq:m10} &&\\
        &S_{im} \geq \sum_{j\in J}s^{\leq}_{jm}\cdot x_{jim} - V\cdot f_{im} && \forall i\in J, m\in M \label{eq:m11} &&\\
        &\sum_{t=0}^{t_{max}}W_{ijmt} = x_{ijm} && \forall i\in J, j\in J, m\in M \label{eq:m12} &&\\
        &\sum_{j\in J}\sum_{t=0}^{t_{max}}t\cdot W_{ijmt}\geq T_{im} && \forall i\in J, m\in M \label{eq:m13} &&\\
        & && \notag &&\\
        &y_{jm}, x_{ijm}, f_{im}, W_{ijmt}\in \{0, 1\} && \forall i\in J, j\in J, m\in M, t = \{0, ..., t_{max}\} \notag &&\\
        &P_{im}, S_{im}, C_{im}, D_{im}, T_{im}\geq 0 && \forall i\in J, m\in M \notag &&
    \end{flalign}
    \endgroup
   
By Constraints (\ref{eq:m1}), each job is assigned to one machine. By (\ref{eq:m2}), if a job is assigned to a machine, then it will be assigned to exactly one slot. Each slot cannot be occupied by more than one job (\ref{eq:m3}). In addition, each machine can only have one firstly-occupied slot (\ref{eq:m4}). Constraints (\ref{eq:m5}) ensure the continuity of each sequence; a slot cannot be occupied if its succeeding one is not already occupied too. Constraints (\ref{eq:m6}) ensure that if a slot $i-1$ is occupied, then the succeeding slot $i$ cannot be the first occupied of machine $m\in M$.

Regarding now the continuous variables, Constraints (\ref{eq:m7}) determine the completion time of all slots, (\ref{eq:m8}) compute the proper deadline value of each slot, (\ref{eq:m9}) define the tardiness of slots, and (\ref{eq:m10})-(\ref{eq:m11}) address to the processing and setup times respectively. Regarding the latter ones, we observe that if a slot $i$ is the first one to be occupied in machine $m$ (i.e., $f_{im} = 1$), then the implied setup time is set to 0 (which is an obvious relaxation for the setup time of the first job $l_{jm}$, for any job $j\in J$)

To compute the value of weighted tardiness, we consider a set of time-indexed binary variables $W_{ijmt}$, indicating that slot $i$ of machine $m$, occupied by job $j$, is delayed by $t$ time units (assuming that $t$ is bounded by a maximum value $t_{max}$). By (\ref{eq:m12}), each triplet $(i,j,m)$ is linked with a value of tardiness $t$, only if the respective assignment variable $x_{ijm}$ is equal with 1. The proper value of $t$ is chosen by (\ref{eq:m13}). The sum of weighted tardiness values is the objective function (\ref{eq:M_obj}).

For completeness of our exposition, we must trivially show the following.

    \begin{lemma}
        $\mathcal{R}$ is a strict relaxation of $\mathcal{P}$. \label{lem:1}
    \end{lemma}
    \begin{proof}
Any feasible solution of $\mathcal{P}$ is a feasible solution for $\mathcal{R}$,
since all constraints of $\mathcal{R}$ are also imposed on $\mathcal{P}$, while the solutions of $\mathcal{P}$ are, in addition, subject to the availability of resources. \end{proof}  

    Although a time-indexed approach is usually not computationally efficient, as even small instances imply a large number of variables and constraints, such formulations are stronger in terms of dual bounds computation. Nevertheless, the tightness of the obtained lower bound is highly instance-dependent: if the setup tasks tend to occupy only a small proportion of the schedule, compared to the processing tasks, then the consumed resources will be available faster and the relaxation of the values of setup times will have a small impact on the minimum objective value. On the contrary, the more sequence-dependent the schedule is, the weaker the lower bound that $\mathcal{R}$ provides will be. 
    
\section{Primal bounds through (meta-)heuristics}
    \subsection{A Constraint Program for resource allocation}
        
Let $\bar{M}$ be a set of given sequences, as provided by variables $x_{ijm}$ in the optimal solution of $\mathcal{R}$. Each sequence $m\in \bar{M}$ has a length $|m|$. Each slot $i\in \{1, 2, ..., |m|\}$ corresponds to an assigned job (i.e., $m_{i}\in J$, $\forall$ $i\in \{1, 2, ..., |m|\}$). If there is no limited resource, a simple re-computation of completion times and tardiness values, replacing the relaxed values of setup times with the (actual) sequence-dependent ones is sufficient. However, if the available resources are not infinite, the following CP model finds an optimal resource allocation.
        \begin{flalign}
            \mathcal{CP}: && \notag &&\\ 
            \text{  min } &\zeta \geq \sum_{j\in J}w_{j}\cdot T^*_{j} && \notag &&\\ 
            &\text{subject to:} && \notag &&\\
            &\texttt{Cumulative}(\texttt{startOf}(\lambda_{m_{i}}), \texttt{sizeOf}(\lambda_{m_{i}}), 1, WR) && \label{eq:s1} &&\\
            &C^*_{m_{i}} = \texttt{endOf}(\lambda_{m_{i}}) + p_{m_{i}m} && \forall m\in \bar{M}, i\in \{1, ..., |m|\} \label{eq:s2} &&\\
            &T^*_{m_{i}} \geq C^*_{m_{i}} - d_{m_{i}} && \forall m\in \bar{M}, i\in \{1, ..., |m|\} \label{eq:s3} &&\\
            &\texttt{startOf}(\lambda_{m_{i}}) \geq C^*_{m_{i-1}} && \forall m\in \bar{M}, i\in \{2, ..., |m|\} \label{eq:s4} &&\\
            &\texttt{sizeOf}(\lambda_{m_{i}}) = s_{m_{i-1}m_{i}m} && \forall m\in \bar{M}, i\in \{2, ..., |m|\} \label{eq:s5} &&\\
            &\texttt{sizeOf}(\lambda_{m_{1}}) = s^{0} && \forall m\in \bar{M} \label{eq:s6} &&\\
            &\lambda_{m_{i}}: \texttt{Interval} && \forall m\in \bar{M}, i\in \{1, ..., |m|\} \notag &&\\
            &C^*_{m_{i}} \geq 0 && \forall m\in \bar{M}, i\in \{1, ..., |m|\} \notag &&\\
            &T^*_{j} \geq 0 && \forall j\in J \notag &&
        \end{flalign}
    
This model employs a set of interval variables $\lambda_{m_{i}}$ for each slot of sequences $m\in \bar{M}$, indicating the start and end time of each setup task. The global \texttt{Cumulative} function ensures that no more than $WR$ setup tasks, which start at $\texttt{startOf}(\lambda_{m_{i}})$ and finish after $\texttt{sizeOf}(\lambda_{m_{i}})$, occupying one resource, will be conducted simultaneously (\ref{eq:s1}). Constraints (\ref{eq:s2}) ensure that the completion time of each job is equal with the termination of the respective setup task interval, increased by the processing time. Constraints (\ref{eq:s3}) determine the value of tardiness. The start time of each setup task should be no less than the completion of the previous job (\ref{eq:s4}). Last, the size of each time interval, by restoring the actual setup time, is defined by (\ref{eq:s5}) and (\ref{eq:s6}).

Exploiting the fact that the presented model computes the optimal resources allocation over a set of sequences in a few seconds, even for instances of considerable scale, we integrate $\mathcal{CP}$ into our implementations of (meta-)heuristics. Notice that solving the MILP relaxation $\mathcal{R}$ and then $\mathcal{CP}$ already provides a feasible solution but of low quality in terms of tardiness (as it is shown in Subsection 4.3).

\subsection{Apparent-Tardiness-Cost-with-Setups (ATCS) Procedure}

The well-known ATCS dispatching rule for the $P| s_{ijm}| \sum_{j}w_{j}T_{j}$ without resource constraints has been introduced in \cite{Pinedo97}. This rule is based on a priority index $I_j$, calculated for each unscheduled job $j$. Here, we use the modified priority index $I_{jm}$, $j \in J$ and $m \in M$ of \cite{Bek19}, after removing the so-called term $CST$, which is the time where a common server is available. We should mention that a similar priority index is also shown in \cite{Chen09} for the $R| s_{ijm}| \sum_{j}T_{j}$ problem. 

Let $J_u$ be the set of unscheduled jobs and $J_s=J\setminus J_u$ be the set of already scheduled jobs, ATCS begins with $J_u = J$ and $J_s = \varnothing$. The first step is to identify the machine $m^{\ast} \in M$ with the minimum load $l_m$ - if more than one machines have equal load, we choose arbitrarily. Then, for each job $j \in J_u$, the following priority index is computed: 

{\scriptsize
\begin{equation}
I_{jm^{\ast}} = \frac{w_j}{p_{jm^{\ast}}} \cdot exp[\frac{-\text{max}(d_j - p_{jm^{\ast}}, 0)}{k_1 \cdot \Bar{p}}] \cdot exp [- \frac{s_{ijm^{\ast}}}{k_2 \cdot \Bar{s}}]
\end{equation}}

\noindent where $\Bar{p}$ and $\Bar{s}$ are the average processing time and setup time respectively. Parameters $k_1$ and $k_2$ are scaling parameters and calculated as: 

{\scriptsize
\begin{alignat}{2}
k_1 = 1.2 \cdot \text{ln}(\mu) - R\\
k_2 = \frac{\tau}{A_2 \cdot \sqrt{\eta}}
\end{alignat}}

\noindent where $\tau$ and $R$ are the due-date tightness and due-date range factors respectively, $\eta = \frac{\Bar{s}}{\Bar{p}}$, $\mu = \frac{|J|}{|M|}$ and $A_2$ is a constant, for which we have $A_2 = 1.8$ if $\tau < 0.8$ or $ A_2 = 2.0$ if $\tau \geq 0.8$. According to \cite{Pinedo97}, two additional modifications need to be made, namely (i) if $\tau < 0.5$, we subtract  $0.5$ from $k_1$ and (ii) if $\eta <0.5 $ and $\mu > 5$, we subtract  $0.5$ from $k_1.$

Job $j^{\ast} \in J_u$ with the maximum priority index $I_{j^{\ast}m^{\ast}}$ is chosen to be scheduled. Similarly to \cite{Bek19}, after choosing job  $j^{\ast}$, we calculate its completion time $C_{ j^{\ast}m}$ per machine $m$, i.e., $C_{j^{\ast}m} = l_m + p_{ j^{\ast}m} + s_{ij^{\ast}m}$, and then schedule $j^{\ast}$ on the machine $m^{\ast\ast}$ imposing the minimum completion time. The weighted tardiness for job $j^{\ast}$ is $wT_{j^{\ast}} = w_{j^{\ast}} \cdot \text{max} (0, l_{m^{\ast\ast}} + p_{j^{\ast}m^{\ast\ast}} + s_{ij^{\ast}m^{\ast\ast}} - d_{j^{\ast}}),$ while the load $l_{m^{\ast\ast}}$ for machine $m^{\ast\ast}$ is increased by  $p_{j^{\ast}m^{\ast\ast}} + s_{ij^{\ast}m^{\ast\ast}}$. Afterwards, job $j^{\ast}$ is moved from $J_u$ to $J_s$.  After scheduling all jobs in this manner, we ran the $\mathcal{CP}$ formulation described in Section 3.1 to deal with the resource constraints. 

\begin{algorithm}[H]
	\fontsize{9}{11}\selectfont
	\caption{\label{alg:atcs} \textsc{ATCS ALGORITHM} }
		\begin{algorithmic}[1]
		\State $ J_u = J, J_s = \varnothing $;
            \State Calculate $k_1, k_2, A_2, \eta , \mu$
		\While {$ J_u \neq \varnothing $}
		\State Find $m^{\ast}$ with \text{min} $l_m$;
		\State Find $j^{\ast}$ with \text{max} $I_{j^{\ast}m^{\ast}}$;
            \State Calculate $C_{j^{\ast}m}$, $ \forall m \in M$;
            \State Assign job $j^{\ast}$ to $m^{\ast\ast}$, where $C_{j^{\ast}m}$ is minimum;
		\State Calculate the weighted tardiness $wT_{j^{\ast}}$ and update $l_{m^{\ast\ast}}$;
  		\EndWhile
		\State Run $\mathcal{CP}$ of subsection 3.1;\\
		\Return $\sum_{j}w_{j}T_{j}$;
	\end{algorithmic}
\end{algorithm}

\subsection{Genetic Algorithm (GA)}
Let us provide some intuition on the design of our GA scheme \cite{Holland1975}, which modifies relevant schemes to accommodate our more complex problem variant. We use the chromosome representation of \cite{Hamz17}, where each gene corresponds to a job. The initial population, i.e., the initial set of solutions, is randomly generated. In order to evaluate the quality of solutions, a decoding algorithm needs to be defined, which transforms the chromosome structure to a feasible job schedule. The decoding algorithm designed is similar to the one in \cite{Hamz17}, but altered due to the most generic nature of our problem. Let $Ch$ be a chromosome, where the decoding algorithm is applied, and $Cut\_off$ be the largest value of $\sum_{j}w_{j}T_{j}$ in the population. Jobs in the chromosome are prioritised, say, from left to right thus, the first job to be scheduled is the one at leftmost. For each job $j$, we calculate the weighted tardiness $w_jT_{jm} = w_j \cdot \text{max} \cdot (0, l_m + p_{jm} + s_{ljm} - d_j)$ on each machine $m \in M$, where $l \in J$ is the last job scheduled on $m$. Job $j$ is assigned on machine $m^{\ast}$, in which the minimum $wT_{jm}$ is achieved. If more than one machines have weighted tardiness values equal to $wT_{jm^*}$, we choose the machine $m$ with the minimum completion time for job $j$: $C_{jm} = l_m + p_{jm} + s_{ijm}$. Then, the load of machine $m^{\ast}$ is increased by $p_{jm^{\ast}} + s_{ijm^{\ast}}$, its last job is set to $j$ and $\sum_{j}w_{j}T_{j}$ is also increased by $w_j \cdot \text{max} \cdot (0, l_{m^{\ast}} + p_{jm^{\ast}} + s_{ijm^{\ast}} - d_j)$. 

Given that we have limited resources, i.e., $WR < |M|$, we must optimally update $\sum_{j}w_{j}T_{j}$ by running our $\mathcal{CP}$ formulation (Subsection 3.1). To accelerate the solution of the algorithm, we impose that $\sum_{j}w_{j}T_{j}$ is lower than $Cut\_off$, otherwise the individual with the chromosome $Ch$ will not survive the next generation. The above condition implies that all solutions, which are not better than the ones already obtained, are excluded from the $\mathcal{CP}$ computation, hence saving a significant amount of computational time. For instances which require the solution of multiple $\mathcal{CP}$ models, the addition of the $Cut\_off$ has also significant impact on obtaining improved solutions. Finally, the value of $\sum_{j}w_{j}T_{j}$ is returned for chromosome $Ch$. The decoding algorithm is presented below.

\begin{algorithm}[H]
	\fontsize{9}{11}\selectfont
	\caption{\textsc{DECODING ALGORITHM} }
		\begin{algorithmic}[1]
		\State Let $Ch$ be a chromosome;
            \State Let $ Cut\_off$ be the maximum weighted tardiness value in the population;
            \State Let $l_m = 0$ be the initial load and $last_m = 0$ the last order on each $m \in M$;
		\For {$ j \in Ch$}
            \For {$m \in M$}
		\State Calculate $w_jT_{jm} = w_j \cdot \max \{0, l_m + p_{jm} + s_{last_{m}jm} - d_j\}$;
            \EndFor
            \State Find $\mathcal{S}_m = \{ m\in M ~|~ wT_{jm}$ is the minimum weighted tardiness$\}$;
            \If {$|\mathcal{S}_{m}|>1$}
            \For {$m \in \mathcal{S}_{m}$}
            \State Choose $m'$ with the minimum $C_{jm'} = l_{m'} + p_{jm'} + s_{last_{m'}jm'}$, where $last_{m'}$ is the last job scheduled in $m'$;
            \EndFor
            \EndIf
            \State Let $m^*$ be the machine chosen from $\mathcal{S}_{m}$;
            \State Assign $j$ to $m^{\ast}$;
            \State Update $l_{m^{\ast}}$, $\sum_{j}w_{j}T_{j}$, $last_{m^{\ast}} = j$;
  		\EndFor
            \If {$WR < |M|$} 
            \If {$\sum_{j}w_{j}T_{j} < Cut\_off$}
		\State Calculate $\sum_{j}w_{j}T_{j}$ by applying $\mathcal{CP}$ of Subsection 3.1;
            \EndIf
            \EndIf\\
		\Return $\sum_{j}w_{j}T_{j}$;
	\end{algorithmic}
  \label{alg:decoding}
\end{algorithm}

The evolution of the population depends on the matching of individuals (parents) inside the population. When two parents are paired, we have to determine how the exchange of information between them (crossover) is performed. The aim of this procedure is to create offsprings, which combine the best parts of the parent chromosomes, thus yielding solutions of better quality. We opt for simplicity and select randomly two individuals from the population as in \cite{Eroglou17}. 
If a randomly generated number $r_c \in [0,1]$ is greater than a predefined parameter $P_c$, the individuals become parents; otherwise, they are returned to the population. Their respective chromosomes are split on a random position (i.e., splitting point). We adopt here the most common and convenient approach of \cite{Eroglou17} by splitting the chromosome in two parts; another possible approach could be a two-point crossover i.e., three parts are created (see e.g., \cite{Hamz17}). Offsprings are created by interchanging the second parts of the parents. It is inevitable that the offspring chromosomes created will be ``faulty", as they may either contain some jobs twice or missing some of them. Thus, we are obliged to ``repair" them by either removing the duplicate jobs at the leftmost and/or by inserting the missing job just after the splitting point. 

\begin{figure}[h!]
    \centering
    \includegraphics[scale=0.35]{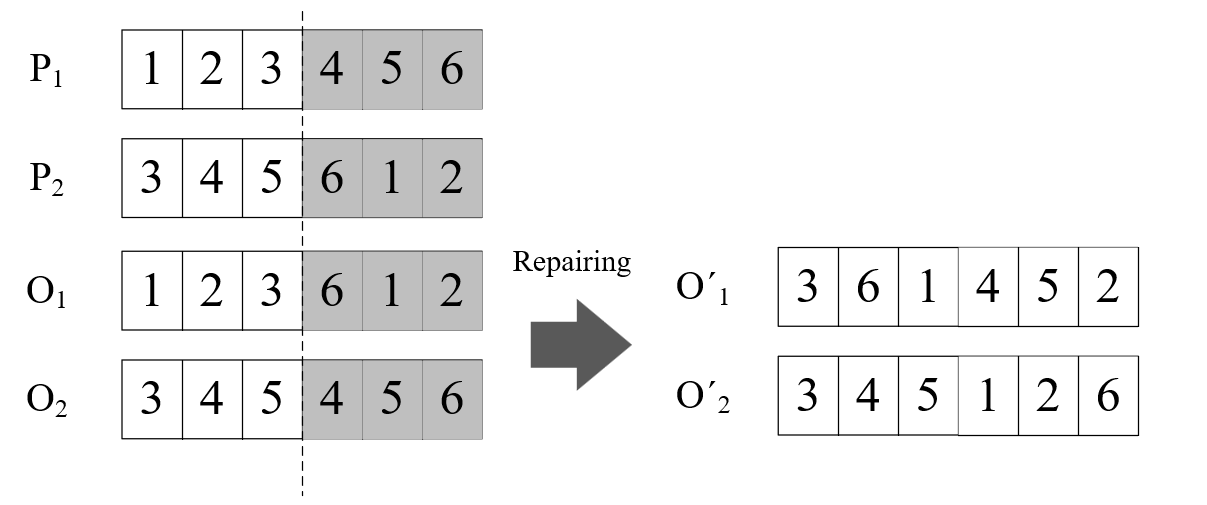}
        \caption{\textbf{Caption:} Crossover operation with chromosome ``repairing"}
    \scriptsize{\textbf{{Figure 1.  Alt Text:}} A diagram which demonstrates the ``crossover" and ``repairing" operations used in the GA}
        \label{fig:ga_schema}
\end{figure}

Figure \ref{fig:ga_schema} shows the crossover operation for two parents $P_1$ and $P_2$. The offsprings $O_1$ and $O_2$ are generated by interchanging the second parts of the parent chromosomes. However, the solutions represented in $O_1$ and $O_2$ are infeasible, since for example, chromosome $O_1$ contains twice genes (jobs) 1 and 2, while not containing genes 4 and 5. Thus, we apply on these chromosomes the repairing procedure, and the final repaired chromosomes are $O^\prime_1$ and $O^\prime_2$. 

\begin{algorithm}[H]
	\fontsize{9}{11}\selectfont
	\caption{\textsc{GENETIC ALGORITHM (GA)} }
		\begin{algorithmic}[1]
            \State Initialise $G_{c} = 0$, where $G_{c}$ is the generation counter;
		\State Generate initial population $P$ of size $P_s$;
            \State Calculate the total weighted tardiness for each individual in $P$ using Algorithm \ref{alg:decoding};
            \While {$ G_{c} \leq G_{max} $}
            \State Create population $P' = P$ and an empty set of new offsprings, $O$;
            \While {$ P' \neq \emptyset $}
            \State Select uniformly at random two individuals $x$ and $y$ from $P'$;
            \State Select uniformly at random $r_c$ from $[0,1]$;
            \If {$r_c > P_c $}:
            \State Use the crossover operation on $x$ and $y$ and generate offsprings $a$ and $b$;
            \State Select uniformly at random $r_m$ from $[0,1]$;
            \If {$r_m < p_m $}:
            \State Use the mutation operation on $a$ and $b$; 
             \EndIf
            \State Calculate total weighted tardiness for $a$ and $b$ using Algorithm \ref{alg:decoding};
            \State Insert $a$ and $b$ to $O$;
            \EndIf
            \State Remove $x$ and $y$ from $P'$;
            \EndWhile
            \State Merge $P$ and $O$ and create new population $P$ with the smallest $|P_s|$ chromosomes in terms of total weighted tardiness values;
            \State Increase $G_c$ by 1;
  		\EndWhile\\
    
		\Return the minimum $\sum_{j}w_{j}T_{j}$ from population $P$;
	\end{algorithmic}
 \label{alg:GA}
\end{algorithm}

In terms of the mutation operation, a random number $r_m \in [0,1]$ is generated and, if $r_m < P_m$, where $P_m$ is the mutation parameter, the chromosome is mutated by randomly swapping two genes (jobs). The fitness value of each offspring is calculated by the Decoding Algorithm (Algorithm \ref{alg:decoding}) and the offsprings are added in the population. Last, the pool of already created chromosomes has to be reduced to a size equal to the initial population. This procedure is called ``survival" and removes chromosomes/solutions with low fitness values, thus driving the population towards better solutions. The algorithm returns the best solution found when a maximum number of generations ($G_{max}$) is reached. 

\subsection{Simulated Annealing (SA)}

The SA algorithm \cite{Metropolis53}, aims at escaping from local optima, by accepting deteriorating neighbour solutions with some probability. This property is the main difference of SA with conventional local search algorithms. SA is influenced by the following factors: an initial temperature $T_0$, a final temperature $T_{cry}$, the cooling rate $q$ and the number of iterations on each temperature $IT$. To run the SA, an initial solution has to be created and be fed to the algorithm. Then, a selected neighboring solution is accepted if its corresponding cost value is smaller than the incumbent one. To avoid being trapped in a local optimum, the selected solution could be also accepted if $p < e^{-\Delta/Tc}$, where $p$ is generated by random uniforms, and valued between 0 and 1. $\Delta$ indicates the difference between the value of the cost function between the current and the new solution and $T_c$ is the current temperature. Thus, when the temperature is slowly reduced by selecting the cooling rate $q$ appropriately, the algorithm converges to the global optimum solution. 

\begin{algorithm}[H]
	\fontsize{9}{11}\selectfont
	\caption{\label{alg:sa} \textsc{SIMULATED ANNEALING (SA)} }
		\begin{algorithmic}[1]
            \State Specify SA parameters $T_0$, $T_{cry}$, $q$, $IT$;
            \State Let $Sol_{cur}$ be the initial solution of SA, obtained from GA or ATCS and $W(T_{cur})$ the total weighted tardiness of the current solution;
            \State Let $Sol_{best}$ be the best solution of the SA algorithm and $W(T_{best})$ its total weighted tardiness;
            \State $Sol_{best} = Sol_{cur}$ ; $W(T_{best}) = W(T_{cur})$ ; $T_c = T_0$, $ext = 1$;
            \While {$ T_{cry} < T_c $}
            \For {$i = 1$ \text{to} $IT$}
            \If {$ ext = 1$}:
            \State Create $Sol_{new}$ from an external swap and set $ext = 0$;
            \State Calculate $W(T_{new})$;
            \Else
            \State Create $Sol_{new}$ from an internal swap and set $ext = 1$;
            \State Calculate $W(T_{new})$;
            \EndIf
            \If {$W(T_{new}) < W(T_{best})$}:
            \State $W(T_{best}) = W(T_{new})$ and $Sol_{best} = Sol_{new}$; 
            \Else
            \State $\Delta = W(T_{new}) - W(T_{cur})$;
            \State Select $U$ uniformly at random from $[0,1]$;
            \State Set $P_{acc} = exp[-\frac{\Delta}{T_c}]$;
            \If {$U < P_{acc}$}:
            \State $W(T_{cur}) = W(T_{new})$ and $Sol_{cur} = Sol_{new}$; 
            \EndIf
            \EndIf
            \EndFor
            \State $T_c = T_c \cdot q$;
  		\EndWhile\\
    
		\Return $W(T_{best})$;
	\end{algorithmic}
 \label{alg:sa}
\end{algorithm}

The initial solution is provided to the SA algorithm from ATCS (Algorithm \ref{alg:atcs}) or GA (Algorithm \ref{alg:GA}). The SA neighbourhood search starts with an external (inter-machine) swap. The corresponding machines are selected based on the value of the cost function - in particular, we choose machines $m'$ and $m''$ with the minimum and the maximum total tardiness values. Then, the most tardy job $j$ in $m''$ is swapped with the least tardy job $i$ in $m'$. Job $j$ is positioned first in $m'$, in order to minimise its delay, while job $i$ is positioned last in $m''$, aiming at not disturbing the current existing sequence. After calculating the new total weighted tardiness value, the algorithm decides whether the new solutions should be accepted or not, based on the criteria described above. It is interesting to note that in Lines 9 and 12, in case where $WR<M$ the calculation of the new total weighted tardiness is performed by applying the $\mathcal{CP}$ of Subsection 3.1.

The next step is to perform an internal (intra-machine) swap between two jobs on the same machine. Firstly, we identify the machine $\hat{m}$ that exhibits the greatest value of total weighted tardiness amongst all machines. Then, we identify jobs $z$ and $k$ with the greatest and smallest value of weighted tardiness, respectively. Finally, we swap $z$ and $k$ on the job sequence of machine $\hat{m}$ and recalculate the new total weighted tardiness. It is important to note that due to the properties and the objective of the problem, it is far from trivial to decide upon which neighbourhood has the greatest prospect of yielding improved solutions. This is why we choose to perform external and internal swaps alternately. Also, this approach ensures that loops of exchanges between the same pairs of jobs will not be performed. 

It should be noted that, this work differs from \cite{Bek19} when exploring neighboring solutions. In \cite{Bek19}, the authors provide six alternative methods of determining a new neighboring solution. More precisely, Methods 1-3 use insertion strategies, while Methods 4-6 use swap operations. Methods 3 and 4 are performed on a completely random manner when selecting jobs and machines. On methods 1-2 and 5-6, even though the machines with the greatest and smallest weighted tardiness values are located, jobs which are swapped or inserted are randomly determined. In addition, new sequences of jobs to machines are also randomly fixed, which may lead jobs with tight deadlines towards the end of sequences, and thus, severely increasing the value of the cost function. Finally, all six methods are examined in every iteration of the algorithm, picking the best one.  

In our approach, we eliminate randomly made decisions and focus only on neighborhoods combining the advantages of our external and internal swaps. Contrary to local improvements (i.e. insertion strategies which may negatively effect the objective cost \cite{Bek19}), our approach aims to drive the algorithm to better solutions, balancing the swapping effects on more than one elements of the schedule. 
Note that, when a job $j$ is removed from a machine $m'$ and inserted to a machine $m''$ (especially if $p_{jm''} > p_{jm'}$) it may not only increase its tardiness (depending on its position on the new sequence), but also the tardiness of jobs succeeding $j$ on $m''$. On the other hand, by swapping jobs we try to avoid machine overloading, while also applying  appropriate machine/job selection strategies.
Let us also note that the idea of selecting the best out of a set of possible methods (as in \cite{Bek19}) would lead Algorithm \ref{alg:sa} to a rapid increase in terms of running time, mainly due to the $\mathcal{CP}$'s multiple calls. These calls cannot be reduced by the cut-off rule used in GA, since the SA algorithm may accept worse solutions by design, and thus, we are obliged to calculate the cost function on every explored neighborhood.

\section{Closing the gap}

All experiments have been performed on a server with 4 Intel(R) Xeon(R) E-2126G @ 3.30GHz processors and 11 GB RAM, running CentOS/Linux 7.0. We used Python 3.8.5 for scripting and Pyomo 5.7.3. The solver of the MILP relaxation $\mathcal{R}$ was CPLEX 20.1 \cite{CPLEX}, via the Pyomo library \cite{Pyomo}, compatible for Python 3.8.5. The $\mathcal{CP}$ model was solved by the CP Optimizer of the DOCplex module, integrated into the provided Python API.

\subsection{Instance Creation}

We generate random instances of $|M| \in \{2,5,10,15,20\}$ machines and $|J| \in \{5 \times |M|, 10 \times |M|\}$ jobs for each value of $|M|$. Since machines are unrelated, the processing time $p_{im}$ of job $i \in J$ in machine $m \in M$ is set to $b_i \cdot a_{im}$, where $b_i$ and $a_{im}$ are selected uniformly at random (u.a.r.) from $[1, 10]$ (as in \cite{Fotakis16}), then increased by a ``noise" selected u.a.r. from $[0, 10]$. For the sequence- and machine-dependent setup times, we set $s_{ijm} = \alpha_{ijm} \cdot p_{jm}$, where $\alpha_{ijm}$ is selected u.a.r. either from $[0.1,0.5]$ or from $[0.5, 1]$ (in a similar fashion to \cite{Lee21}). We also use an alternate manner to generate setup times as $U(5,25)$, where $U$ is the uniform distribution, as in \cite{Bek19}. Weights, denoted by $w_j$, are generated u.a.r from $[1, 10]$.

Regarding the values of deadlines $d_{j}$, we have to determine the due date tightness factor $\tau$ and the due date range factor $R$, so as to create the corresponding value $d_{j} = U(C_{max}\cdot (1 - \tau - \frac{R}{2}), C_{max}\cdot (1 - \tau + \frac{R}{2}))$. The makespan is calculated as $C_{max} = \frac{\sum_{j}({P^1_{j}})+ S^1}{|M|} $, where $P^1_{j} = \text{min}_{m\in M}\{p_{jm}\}$. The computation of $S^1$ takes advantage of the zero initial setup time. Firstly, we define $S_j = \text{min}_{i\in J\setminus \{i\}, m\in M}\{s_{ijm}\}$. These values are sorted in ascending order and, since no setup time is required for the first $|M|$ jobs,  we assume that $S^1$ is equal to the sum of the $S_j$ for the first $|J|-|M|$ jobs \cite{Bek19}. Our instances are generated for $\tau = \{0.5, 0.8\}$ and $R = 0.8$.

 Each instance is solved for $WR \in \{\ceil{\frac{|M|}{2}} , |M| \}$, where $WR$ is the number of working resources. Each (meta-)heuristic algorithm is solved five times, in order to remove stochasticity of solutions. The above yield a total number of 600 experiments, which when multiplied by the 4 different methods applied rises to 2400 distinct executions or, given our time limits, more than $500$ hours. 

\subsection{Implementation details}
The performance of the MILP lower-bounding relaxation is instance-specific. That is, if a high variance between the relaxed and actual values of setup times is noticed, the optimal objective value of $\mathcal{R}$ will be significantly distant from the respective primal bound. In addition, due to the set of four-indices variables $W_{ijmt}$ induces, the MILP relaxation itself cannot reach optimality for instances of more than 50 jobs within one hour (doubling that time limit creates no significant size increase). Even this can be accomplished only by using the best primal bound obtained by our (meta-)heuristics work as warm-starts for $\mathcal{R}$ \cite{CPLEX}. 

We also obtain the respective primal solution implied by $\mathcal{R}$. If $WR = |M|$, then we compute an upper bound by simply replacing the relaxed values of setup times with the actual ones for the fixed sequences (denoted by $\bar{M}$). If $WR < |M|$, then we compute a solution by applying our $\mathcal{CP}$ which requires no more than very few seconds to be solved to optimality. The procedure of obtaining primal solutions from sequences which have been constructed by $\mathcal{R}$ is simply referred as $MIP$ algorithm hereinafter.

The GA parameters i.e. Population Size ($P_s$) and Number of Generations ($G_{max}$) were set to 100 and 150 respectively. Since the SA algorithm is computationally expensive, due to the number of $\mathcal{CP}$s to be solved, its parameters are set to: $T_0 = 500$, $T_{cry} = 1$, $q = 0.9$ and $IT = 50$. Note that the values for the GA and SA parameters are close to the values used in \cite{Hamz17}, yet carefully reduced to remain within the boundaries of reasonable computational times. For the same reason, we impose a time limit of 3600 seconds for the GA algorithm. In practice, this means that some instances will run for a smaller number of generations, however, the population size $P_s$ remains the same.

\subsection{Experiments on benchmark instances}
Table \ref{tab:MH_res} shows the results. Column `\emph{MIP}' shows the sequential solution of the MILP relaxation $\mathcal{R}$ and the $\mathcal{CP}$ model, where  $\mathcal{R}$ is computed for instances of up to $50$ jobs and $10$ machines within the 1-hour time limit. Regarding the implementation of the (meta-)heuristics, the presented values are averages of 5 consecutive runs on every instance and the table aims at comparing the results of the proposed methods under each possible combination of size of instances ($J\times M$), job to machine dependency ($|M|$), setup times ($\alpha$), setup resources ($WR$) and different tightness of deadlines ($\tau$). `Gap' is calculated as: $\textit{Gap} =\frac{\textsc{Alg\_sol - Lower Bound}}{\textsc{Alg\_sol}}$ and it is provided exclusively for the instances that admit a solution of $\mathcal{R}$. \textit{Time} is the computational time (in seconds). Finally, $Err$ is equal to $\frac{\textsc{Alg\_sol - Alg\_best\_sol}}{\textsc{Alg\_sol}}$, where $Alg\_sol$ is the solution of an algorithm and $Alg\_best\_sol$ is the best solution found by any of the (meta-)heuristics.   

\begin{table} [H]
    \fontsize{6}{12}\selectfont
    \tbl{Algorithm results on benchmark instances}
    {\resizebox{\textwidth}{!}{
        \begin{tabular}{|c|c|c|c|c|c|c|c|c|c|c|c|c|c|c|c|c|} \hline 
            \multicolumn{2}{|c|}{\multirow{3}{*}{Instance}} &                      
              \multicolumn{15}{|c|}{Method} \\ 
              \cline{3-17}
            \multicolumn{2}{|c|}{}   & \multicolumn{3}{|c|}{MIP}& \multicolumn{3}{|c|}{GA} & \multicolumn{3}{|c|}{GA+SA} & \multicolumn{3}{|c|}{ATCS}  & \multicolumn{3}{|c|}{ATCS+SA} \\ \cline{3-17}
            \multicolumn{2}{|c|}{}  & Gap & Time & Err &Gap & Time &  Err & Gap & Time & Err & Gap& Time& Err & Gap & Time & Err \\ \cline{1-17}
            \multirow{10}{*}{$J\times M$} & 10x2& 19.53 & 34 & 6.4 & 17.47 & 206  & 4.74 &  17.18 & 113 & 4.44  & 40.69 & $<$ 1 & 31.99 & 38.53 & 123 & 29.22\\
            & 20x2 & 51.36	&3600	&22.38&	40.03	&1195	&4.19	&40.02&	967	&4.17	&51.17&	$<$ 1	&21.55&	50.14&	857	&19.85\\
            & 25x5 & 69.24	 &3600	 &14.53 &	66.44 &	327 &	8.16 &	66.36 &	138	 &7.87 &	74.02 &	$<$ 1	 &29.51 &	73.88 &	247	 &29.01 \\
            & 50x5 & 91.1	&3600	&	40.24&	86.59	&1193	&10.2	&86.44	&	934	&8.88	&87.03	& $<$ 1	&	12.28&	86.87	&	847	&10.66 \\
            & 50x10 & 79.89 &	3600 &	5.73&	79.7&	1229&	10.15&	79.6 &	949	& 9.55 & 83.51& $<$ 1	&24.55 & 83.37 &1062 &24 \\
            & 100x10 & - &  - & - & - &  1534  & 28  &- & 1276 & 26.22 & -& $<$ 1 &4.6 &- &1275 & 4.11 \\
            & 75x15 & - & -  & - & - &  1182  & 12.08  & -& 735 & 11.16 & -& $<$ 1 &19.09 & -&933 & 18.66 \\
            & 150x15 & - & -  & - & - &  1704  & 46.09  &- & 1525 & 44.61 &- & $<$ 1 &4.47 &- &1507 & 4.36 \\
            & 100x20 & - & -  & - & - &  1365  & 19.75  &- & 868 & 18.55 & -& $<$ 1 &13.35 &- & 1271 & 12.36\\
            & 200x20  & - &  - & - & - &  1737  & 58.62  & -& 1691 & 57.9 &- & $<$ 1 & 5.71 &- & 1794 & 5.33\\
            \hline
            \multirow{2}{*}{$M$}  & $5 \cdot |M|$ &  56.22 &2411 &8.88 &54.54 &862 &10.97 &54.38 &560 &10.31 &66.07 & $<$	1 &23.7 &65.26 &727 &22.65 \\
            &  $ 10 \cdot |M|$   & 71.23 &3600&	31.31&	63.31&	1473&	29.42&	63.23& 1271& 28.36&	69.1& $<$ 1  &9.72 &68.5&	1256&	8.86 \\
             \hline
            \multirow{3}{*}{$\alpha$}& [0.1,0.5]  & 59.02&2885&	11.96& 56.55&	637&19.23&56.47&407&18.51&66.17&$<$	1&15.07&65.85	&353	&14.73\\
             & [0.5,1]  & 61.01& 2890 &21.83&56.56&	1496&	19.64&	56.52&	1268&	18.63&	67.24& $<$	1	&23.45&	65.83&	1351&	21.55 \\
            & U (5,25)  & 66.65	&2889&	19.77&61.03&1369& 21.72 &60.76 &1073& 20.87 &68.44 & $<$	1 &11.62 &68 &1271	&11\\
            \hline
            \multirow{2}{*}{$WR$}  & $\ceil{\frac{|M|}{2}}$ &  63.73 &2887& 17.27 &61.08 &2317 &21.06 &61.01 &1830 &20.21 &70.09&	$<$ 1 &17.34 &68.85 &1980 &15.79 \\
            &  $|M|$   & 60.72 &2887 &	18.44	&55.02&	18	&19.34&	54.83&	2&	18.46	&64.48& $<$	1	&16.08	&64.26&	2	&15.73\\
             \hline
             \multirow{2}{*}{$\tau$} & 0.5  & 67.09 &2884 &18.35&63.45 &1157&	29.95 &63.35 &911 &28.79 &74.1 & $<$ 1 &18.55 &72.72 &1008 &16.96\\
             & 0.8 &  57.36 &2890 &17.36 &52.64 &1178 &10.45 &52.48 &921 &9.88	& 60.47 & $<$ 1 &14.87 &60.4 &975&14.55 \\
            \hline
            \multicolumn{2}{|c|}{Mean} & 62.22 & 2887   & 17.86 & 58.05 & 1167 & 20.2 & 57.92  & 916 & 19.33 & 67.29 & $<$ 1 & 16.71 & 66.56 & 992 & 15.76 \\
            \hline
        \end{tabular}
    }}
    \label{tab:MH_res}
\end{table}

The values of Gap appear relatively high, even on instances of $20\times 2$ for all the methods. This does not necessarily correspond to the quality of solutions, as it could be attributed also to the weakness of the lower bounds found by $\mathcal{R}$. On the other hand, as already stated above, it remains important to have even an estimation for gaps, especially when taking under consideration not only the difficulty of the problem handled here for the first time but also the considerable size of our instances. Another critical observation is that the SA, whether receiving the solution from GA or ATCS, can only slightly improve it, within an average of 950 seconds. It should, also be noted that since the SA solutions depend on the solutions of the other algorithms, the time needed on each execution is the sum of finding a solution and then using SA. For example, GA+SA solutions for the 10$\times$ 2 instances are provided in $206 + 113 = 319$ seconds, and, as the instance size grows, so does the total time. 

The most insightful outcome shown by Table \ref{tab:MH_res} concerns the comparison of GA to ATCS. It is obvious that for instances up to $75$ jobs  and $15$ machines, GA performs better in terms of both $Gap$ and $Err$. However, when instances grow beyond $100$ jobs, ATCS severely outperforms GA. An explanation could be that the larger instance size yields an exponentially larger solution set, which a meta-heuristic cannot explore effectively, even with the advantage of using more computational time. 

Let us mention two further facts for the GA algorithm. Firstly, for instances having $WR < M $ thus requiring solving $\mathcal{CP}$ for each solution, the average number of generations run is 134, meaning that for several such instances the GA timed out due to its 1-hour limit without being able to reach the $G_{max}=150$ generations. In addition, on average 6743 different solutions were examined from the GA algorithm on each experiment. The rule proposed in Algorithm \ref{alg:decoding}, which skips the $\mathcal{CP}$ execution based on the $Cut\_off$ was activated 2037 times on average (approximately 30 \%). When considering that the duration of each $\mathcal{CP}$ iteration is almost 1 second, it is understood that this cutting-off rule is valuable for GA. The percentage of computational burden coming from $CP$ component is also revealed when comparing the corresponding times between instances of $WR = |M|$ and  $WR = \ceil{\frac{|M|}{2}} $. As shown in Table \ref{tab:MH_res}, all algorithms (excluding $MIP$) are finished under a few seconds on average when $WR = |M|$, while this time grows immensely -up to 1 hour on some instances- when $WR = \ceil{\frac{|M|}{2}}$. In short, resource constraints are what makes instances really hard for GA, as $\mathcal{CP}$ becomes a prerequisite. This is in sharp contrast to the variant of \cite{Bek19} with a common server, which is the closest variant to ours in the current literature.

Furthermore, GA performs better than ATCS on instances with the following characteristics; $|J| = 5 \cdot |M|$, $ \alpha \in [0.5, 1]$, which implies that setup times are close to the values of the processing times and $\tau = 0.8$, which corresponds to tighter deadlines. This is quite interesting as the instances $\alpha \in [0.5, 1]$ and $\tau = 0.8$ seemed to be the most challenging ones, significantly tightening the constraints of the problem. Another interesting fact is that, for instances where $J\times M = 50\times 10$, MIP exhibits the lower value of $Err$ with 5.73 \%, more than 60\% better than the second best method. However, it should be clarified that MIP should not be directly compared to the other algorithms (apart from $J \times M$ rows), since the average values provided are calculated over only a limited size of instances (up to $J\times M = 50\times 10$). This, also, explains why the overall average values of MIP are lower than the ones showed for the rest of the algorithms.
\subsection{Experiments on real data}
The ultimate goal of designing efficient scheduling algorithms is to apply them on problems faced by actual manufacturers. Since 
 it is inevitable to observe differences on the composition of benchmark and real data, we evaluate the performance of our computational tools on a real industrial setting, originating from the textile industry. In particular, we are motivated from the case of textile manufacturer PIACENZA, located in Northern Italy. The industrial setting of PIACENZA consists of two types of parallel looms, which are assigned to produce orders so that the required demand is covered. To be consistent with the general nomenclature of scheduling problems, we refer to \emph{looms} and \emph{orders} as machines and jobs respectively, hereinafter. The remaining aspects of the problem are as described in Section \ref{Section:Problem}; each job is featured with a weight value $w_{j}$ and a deadline $d_{j}$. Specifically for the $WR$ parameter, the plant operates with human workers who are responsible for setting up the machines so as to produce the amount of demanded fabric. 

The textile manufacturer has provided us with a set of 45 weekly instances, spanning from January to December 2020. The size of the real instances is of comparable scale with the randomly generated benchmarks, as the number of jobs ranges from $7$ to $94$. Currently $WR=3$, while the number of machines $|M|$  ranges from $5$ to $10$. Next, we  compare our proposed algorithms in order to decide upon the most suitable approach for PIACENZA.

Table \ref{tab:real_res} presents the results on historical instances. It should be noted that, since the SA algorithm demands a large amount of computational time without delivering considerable improvements, it is excluded from the experimentation on real instances. In addition, due to the presence of several extremely tight deadlines, the time-indexed formulation of $\mathcal{R}$ is not applicable, even for the smallest instances (i.e., a large number of $t$ distinct time instances would yield the construction of a respective number of variables $W_{ijmt}$). Consequently, no \emph{Gap} values are obtained. We also note that $\tau$ is defined as $\tau = 1 - \frac{\Bar{d}}{C_{max}}$, with $\Bar{d}$ being the average due date and $C_{max}$ being the estimated makespan. The due date range factor ($R$) is defined as $R = \frac{d_{max} -d_{min}}{C_{max}}$, where $d_{max}$ and $d_{min}$ are the maximum and minimum deadlines respectively.

  \begin{table}[h!]
  \tbl{Algorithm results on real instances} 
  {\begin{adjustbox}{width = \textwidth}
 \begin{tabular}{|c|c|c|c|c|c|c|c|c|c|c|c|c|c|c|c|c|c|c|c|c|} \hline
 \multirow{2}{*}{\begin{tabular}{c}\\ $|J|$\end{tabular}} 
 &\multicolumn{2}{c|}{\% Err GA}  &  \multirow{2}{*}{\begin{tabular}{c} Mean \\t(s)\end{tabular}}  & \multicolumn{2}{c|}{\% Err ATCS}  &  \multirow{2}{*}{\begin{tabular}{c} Mean \\t(s)\end{tabular}} & \multirow{2}{*}{\begin{tabular}{c}\\ $|J|$\end{tabular}} 
 &\multicolumn{2}{c|}{\% Err GA}  &  \multirow{2}{*}{\begin{tabular}{c} Mean \\t(s)\end{tabular}}  & \multicolumn{2}{c|}{\% Err ATCS}  &  \multirow{2}{*}{\begin{tabular}{c} Mean \\t(s)\end{tabular}}\\ \cline{2-3}\cline{5-6} \cline{9-10}\cline{12-13} 
 & 5& 10& & 5 & 10 & & & 5& 10 &  & 5 & 10 &  \\ \hline
7 & 0 & 0 &430 &  0.26 & 0 & $<$ 1 & 50 & 3.89 & 0.03 & 990 & 59.58 & 36.77 & $<$ 1 \\ 
16 & 0 & 0 & 78 & 0 & 0 & $<$ 1 & 50 & 11.43 & 0 & 706 & 43.33 & 32.16 & $<$ 1   \\ 
17 & 0.23 & 0.04 & 66&  0.26 & 0.51 & 1 &  51& 1.68 & 0.29 & 84 & 0.09 & 0.11 & 1 \\ 
21 & 0.01 & 0 & 498 & 1.11 & 0.5 & $<$ 1 &  51& 2.05 & 0 & 1208 & 25.52 & 13.54& 1 \\ 
26 & 1.25 & 0.07 & 20 & 29 & 16.84 & $<$ 1 & 51 & 3.33 & 4.68 & 163& 17.79 & 19.35 & 1 \\ 
26 & 0 & 0 & 527 & 97.04 & 100 & $<$ 1 & 56 &  36.99& 0 & 313& 79.77 & 80.9& 1 \\ 
26 & 24.17 & 8.7 & 26 & 60.22 & 36.24 & $<$ 1 &  57& 3.37 & 0.16 &1801 & 41.94 &21.38 & 1 \\ 
29 & 0.6 & 0.07 & 18 & 29 & 16.84 & $<$ 1 &  57& 0.1 & 0.01 & 1261 & 15.72 & 5.55& 1 \\ 
32 & 0.61 & 0.16 & 44 & 9.75 & 3.76 & $<$ 1 &  58& 16.72 & 10.13 & 1274 & 46.67 & 56.58& 1 \\ 
35 & 1.9 & 0 & 81 & 57.97 & 23.78 & $<$ 1 & 58 & 1.63 & 1.31 & 3600 & 15.81 &24.03 & 1 \\ 
37 & 0 & 0 & 39 & 56.28 & 21.72 & $<$ 1 &  60&  0.03& 0 & 931& 35.52 & 15.29& $<$ 1 \\ 
38 & 4.22 & 0 &347 & 62.41 & 28.2 & $<$ 1 &  60& 0 & 0 & 801 & 100 & 100 & 1 \\ 
40 & 21.78 & 0 & 498 & 70.23 & 39.13 & $<$ 1 & 61 & 4.45 & 0 & 561 & 25.05 &8.86  & $<$ 1 \\ 
40 & 3.24 & 0.53 & 123 & 0.18 & 0.23 & 1 &  65& 5.58 & 2.88 & 3600 & 14.98 & 26.11 & 1 \\ 
43 & 1.48 & 0.33 & 32 & 1.41 & 2.06 & 1 &  65&  8.47& 5.4 & 1478& 31.54 & 33.67& 1 \\ 
44 & 1.2 & 0 & 632 & 28.08 & 13.51 & $<$ 1 &  67& 6.28 & 0 & 1033 & 83.23 & 73.25 & 1 \\ 
44 & 0.34 & 0.14 & 49 & 1.51 & 0.92 & 1 &  67&  6.87 & 3.61  & 3600& 14.81  & 23.48 & 1 \\ 
44 & 0.04 & 0 & 824 & 1.92 & 0.73 & 1 &  69&  3.03 & 0.45 & 2676& 10.23 & 6.68 & 1 \\ 
45 & 0.43 & 0.11 & 27 & 1.32 & 0.64 & $<$ 1 &  72& 0.63 & 0 & 1442 & 22.86 & 6.05 & 1 \\ 
46 & 0.26 & 0.26 & 3600 & 2.74 & 2.56 & 1 &  74& 26.38 & 6.2 & 201 & 2.4 & 8.9 & 1 \\ 
46 & 0.97 & 0.17 & 65 & 0.38 & 0.49 & 1 &  78&  8.76 & 2.02 & 261 & 18.48 & 8.9 & 1 \\ 
48 & 3.89 & 0.03 & 990 & 59.58 & 36.77 & $<$ 1 & 94 & 10.79 & 4.23 & 3151 & 22.53 & 46.86 & 1  \\ 
49 & 0.96 & 0.03 & 1173 & 45.03 & 28.71 & 1 &  -& - & - & -& - &- & - \\ 
\hline
Mean &  2.94 & 0.46 & 444 & 25.51  & 15.54 & $<$ 1 & - & 7.38 & 1.91& 1379  & 30.4&  27.88 & 1  \\ 
\hline
\end{tabular}
\end{adjustbox}}
  \label{tab:real_res}
\end{table}

From Table \ref{tab:real_res}, it is obvious that the GA algorithm returns better solutions than ATCS on all instances, in sharp contrast with benchmark instances of similar size. More precisely, the $ \% Err$ of GA is 5.11 \% on average for 5 machines and 1.17 \% for 10 machines, while the $ \% Err$ of the ATCS algorithm is 27.9 \% and 21.57 \% respectively. Timewise, as already discussed for the benchmark experiments, the GA algorithm consumes more computational time on average (901 seconds) compared to the ATCS (1 second). For the smaller size instances of the real data ($|J| < 50$), the $\% Err$ of GA is 1.7 \% on average, while the corresponding value for the ATCS algorithm is 20.53 \%. We obtain similar results, when $|J| \geq 50$. The $\% Err$ of GA is 4.65 \% on average, while ATCS returns solutions of 29.14 \% in terms of $\% Err$. 

Interestingly, experimentation on real data provides quite different results than the ones with the benchmark data. More precisely, the average $\% Err$ of GA on real data for instances of 25x5, 50x5 and 50x10 (4.65\%) is less than half of the $\% Err$ in benchmark experiments of the same size (9.5\% approximately). The reverse pattern is exhibited by the ATCS algorithm, as the average $\% Err$ on real data (29.14 \%) is greater than the one on benchmark data (22.11 \% approximately). 

This observation can safely be associated with the fact that, on the real data, the tardiness parameters $\tau$ and $R$ obtain values drastically different than the ones generated u.a.r. from $[0, 1]$: $\tau=-74$ and $R=909$. As a result, the scaling parameters $k_1$ and $k_2$ of the ATCS algorithm, obtain extreme values and subsequently, the priority index is calculated in a radically different manner. On the positive side, we may deduce that ATCS and GA perform in a complementary manner, hence forming together an encouraging and versatile approach. 
\section{Concluding insights}
The importance of the weighted tardiness objective lies on its connection to customer satisfaction. The cost of losing a client to competition might be unbearable even for successful businesses and tardiness could determine such an outcome. This applies, indicatively, to the highly competitive sector of textile manufacturing. As a result, managers need innovative computational tools to ensure that production management in the plant is performed at the highest level of efficiency. These facts underline that the methods proposed should be examined in multiple terms, namely for (i) assigning the incoming orders to machines in a dynamic manner and (ii) serving as a fast simulation tool for the plant managers, rendering them able to alter the environment settings, e.g., the number of resources. These capabilities facilitates both operational decision making process and strategic foresight on investment opportunities.  
 
Indicatively, the PIACENZA team has examined a possible increase on the number of machines $|M| = 15$ and also, an increase on the number of working resources $WR = 5$. This will simply estimate the potential profit, in terms of tardiness minimisation, from integrating more machines or more workers into the plant floor. Under the textile manufacturing context, this translates into buying or even leasing looms under short term contracts or/and adding more workers to the plant. 

This additional extensive experimentation on real data has provided some insightful results. Increasing the number of workers, i.e. changing $WR$ from 3 to 5, decreases weighted tardiness by approximately 22 \%. On the other hand, an increase on the number of machines $|M|$ from 10 to 15 decreases the weighted tardiness only by 3.2 \%. Thus, it seems that an investment on new machines is not justified, while hiring new workers may be in the best interest of the plant. However, we should note that by studying the data provided, an additional important aspect to focus on is the assignment of appropriate due dates on each incoming order \cite{chen2022}, since approximately the 15 \% of the incoming orders arrives already late, i.e. past its due date to start its processing, while 17 \% has been assigned with a due date greater than 30 days.

Since this is the first study on weighted tardiness under resource constraints for sequence-dependent setups on unrelated machines, we believe it could motivate further work in several aspects.  From the technical viewpoint, we underline the importance of improving lower bounds. This might translate into designing better and more efficient MILP formulations for a (less) relaxed version of the initial problem or in exploiting some currently undermined properties of the problem. It is encouraging that our lower-bounding MILP offers worst-case gaps on instances of size up to 50 jobs and 10 machines for the weighted tardiness objective. 

An exact approach would be worthwhile to pursue, even for medium-sized instances. At the same time, as we aim at methods and algorithms dealing with the continuously increasing requirements of modern manufacturing, the extension of the size of instances which we are able to solve is of great importance. 

From a managerial viewpoint, it is vital to support the decision making process in additional plants, also because our experimentation has already shown different computational evidence when shifting from benchmark to real instances. Specifically, we should examine weighted-tardiness problems in make-to-order environments that operated under tight deadlines from a three-sided perspective: (i) the efficient planning and scheduling of new incoming orders, (iib) the design of algorithms and methods to be utilised as simulation tools by diverse company personnel, and (iii) the extensive experimentation in order to provide more in-depth strategic insights in terms of managing machines, resources and any other constraint imposed or dictated by the plant. 


\section*{Acknowledgement}

This research has been supported by the EU through the FACTLOG Horizon 2020
project, grant number 869951.

\section*{Conflict of interest statement}

The authors declare that they have no known competing financial interests or personal relationships that could have appeared to influence the work reported in this paper.

\section*{Data availability statement}

Instances, results and codes that support the findings of this work are available at \url{https://github.com/svatikiot/Weighted_Tardiness_Experiments}.

\setstretch{1.25}


\begin{thebibliography}{8}
    \bibitem{IJPR}
    Avgerinos, I., Mourtos, I., Vatikiotis, S., Zois, G.,
    \newblock Scheduling unrelated machines with job splitting, setup resources and sequence dependency.
    \newblock {\em International Journal of Production Research}, (2022)

    \bibitem{MIM}
    Avgerinos, I., Mourtos, I., Vatikiotis, S., Zois, G.,
    \newblock Exact methods for tardiness objectives in production scheduling.
    \newblock {\em IFAC-PapersOnLine, 55:10}, 2487-2492, (2022)
    
    \bibitem{Bek19}
    Bektur, G., Saraç, T.
    \newblock A mathematical model and heuristic algorithms for an unrelated parallel machine scheduling problem with sequence-dependent setup times, machine eligibility restrictions and a common server.
    \newblock {\em Computers \& Operations Research, 103}, 46-63, (2019).
    
    \bibitem{Bilyk}
    Bilyk, A., Mönch, L.,
    \newblock A Variable Neighborhood Search Approach for Planning and Scheduling of Jobs on Unrelated Parallel Machines.
    \newblock {\em Journal of Intelligent Manufacturing, 23}, 1621–1635 (2012).

    \bibitem{Chen09}
    Chen, J.F.,
    \newblock Scheduling on unrelated parallel machines with sequence-and-machine dependent setup times and due-date constraints.
    \newblock {\em The International Journal of Advanced Manufacturing Technology, 44}, 1204-1212, (2009).

    \bibitem{chen2022}
    Chen, R., Gao Y., Geng, Z., Yuan, J.
    \newblock Revisit the scheduling problem with assignable or generalized due dates to minimize total weighted late work.
    \newblock {\em International Journal of Production Research}, 1-19, (2022).
    
    \bibitem{Emmons87}
    Emmons, H.,
    \newblock  Scheduling to a common due date on parallel uniform processors. 
    \newblock {\em Naval Research Logistics, 34:6}, 803-810 (1987).

    \bibitem{Eroglou17}
    Eroglu, D.Y.,  Ozmutlu, H.C.,
    \newblock  Solution method for a large-scale loom scheduling problem with machine eligibility and splitting property.
    \newblock {\em The Journal of the Textile Institute, 108:12}, 2154-2165 (2017).

    \bibitem{Fotakis16}
    Fotakis, D., Milis, I., Papadigenopoulos, O., Vassalos, V., Zois, G.,
    \newblock Scheduling MapReduce jobs on identical and unrelated processors.
    \newblock {\em Theory of Computing Systems, 64:5}, 754-782 (2016).	
    
    \bibitem{Grah79}
    Graham, R.L., Lawler, E.l., Lenstra, J.K., Kan, A.H.G.R.
    \newblock Optimization and Approximation in Deterministic Sequencing and Scheduling: a Survey.
    \newblock {\em Annals of Discrete Mathematics, 5}, 287-326, (1979).

    \bibitem{Hamz17}
    Hamzadayi, A., Yildiz, G.
    \newblock  Modeling and solving static m identical parallel machines scheduling problem with a common server and sequence dependent setup times
    \newblock {\em Computers and Industrial Engineering, 106}, 287-298 (2017).

    \bibitem{Holland1975}
    Holland, J. H.
    \newblock Adaptation in natural and artificial systems
    \newblock {\em Michigan: The University of Michigan Press}, (1975).

    \bibitem{Hook07}
    Hooker, J.N.,
    \newblock Planning and Scheduling by Logic-Based Benders Decomposition.
    \newblock {\em Operations Research, 55:3}, 588-602 (2007).


    \bibitem{CPLEX}
    IBM ILOG CPLEX Optimization Studio CPLEX User’s Manual, Version 12 Release 8,
    \newblock Starting from a solution: MIP starts.
    \newblock {\em \url{https://www.ibm.com/docs/en/SSSA5P_12.8.0/ilog.odms.studio.help/pdf/usrcplex.pdf}}, 255-260.

    \bibitem{Iori22}	
    Iori, M., Locatelli, A., Locatelli, M.
    \newblock A GRASP for a real-world scheduling problem with unrelated parallel print machines and sequence-dependent setup times.
    \newblock {\em International Journal of Production Research}, (2022).

    \bibitem{Jiang22}
    Jiang, Z., Yuan, S., Ma, Jing, Wang, Q.,
    \newblock The evolution of production scheduling from Industry 3.0 through Industry 4.0.
    \newblock {\em International Journal of Production Research, 60:11},  3534–3554, (2022).

    \bibitem{Kim20b}
    Kim, J, Kim, H. J.,
    \newblock Parallel machine scheduling with multiple processing alternatives and sequence-dependent setup times
    \newblock {\em International Journal of Production Research, 59:18},  5438-5453 (2020).

    \bibitem{Lee21}
    Kim, H.J., Lee, J.H.,
    \newblock Scheduling uniform parallel dedicated machines with job splitting, sequence-dependent setup times, and multiple servers.
    \newblock {\em Computers \& Operations Research, 126}, (2021).

    \bibitem{Kim20a}
    Kim, J.G., Song, S., Jeong, B.
    \newblock Minimising total tardiness for the identical parallel machine scheduling problem with splitting jobs and sequence-dependent setup times.
    \newblock {\em International Journal of Production Research, 58:6}, 1628-1643, (2020).

    \bibitem{Kuh16}
    Kuhpfahl, J., Bierwirth, C.,
    \newblock A study on local search neighborhoods for the job shop scheduling problem with total weighted tardiness objective.
    \newblock {\em Computers \& Operations Research, 66}, 44-57, (2016).

    \bibitem{Lawler93}	
    Lawle,r E.L., Lenstra, J.K., Rinnoy Kan, A.H.G., Shmoys, D.B.,
    \newblock Chapter 9: Sequencing and scheduling: Algorithms and complexity.
    \newblock {\em Handbooks in Operations Research and Management Science, 4}, 445-522, (1993).

    \bibitem{Lee18}
    Lee, C.H.,
    \newblock A dispatching rule and a random iterated greedy metaheuristic for identical parallel machine scheduling to minimize total tardiness.
    \newblock {\em International Journal of Production Research, 56:6}, 2292–2308, (2018).

    \bibitem{Lee13}
    Lee, J.H., Yu, J.M., Lee, D.H.,
    \newblock A tabu search algorithm for unrelated parallel machine scheduling with sequence and machine-dependent setups: minimizing total tardiness.
    \newblock {\em The International Journal of Advanced Manufacturing Technology, 69}, 2081-2089, (2013).

    \bibitem{Pinedo97}
    Lee, Y.H., Pinedo, M.,
    \newblock Scheduling jobs on parallel machines with sequence-dependent setup times.
    \newblock {\em European Journal of Operational Research, 100}, 464-474, (1997).
    
	\bibitem{Len77}	
    Lenstra, J.K., Rinnooy Kan, A.H.G., Brucker, P.
    \newblock Complexity of Machine Scheduling Problems.
    \newblock {\em Annals of Discrete Mathematics, 1}, 343-362, (1977).
    
    \bibitem{Liao22}	
    Liao, B., Lu, S., Jiang, T., Zhu, X.,
    \newblock A variable neighborhood search and mixed-integer programming models for adistributed maintenance service network scheduling problem.
    \newblock {\em International Journal of Production Research}, (2022).

    \bibitem{Liaw03}	
    Liaw, C.F., Lin Y.K., Cheng, C.Y., Chen, M.
    \newblock Scheduling unrelated parallel machines to minimize total weighted tardiness.
    \newblock {\em Computers \& Operations Research, 30}, 1777-1789, (2003).

    \bibitem{Metropolis53}	
    Metropolis, N., Rosenbluth, A. W., Teller, A. H., Teller, E.
    \newblock Equation of state calculation by east computing machines.
    \newblock {\em Journal of Chemical Physics, 21}, 1087–1091, (1953).

    \bibitem{Mos21}	
    Moser, M., Musliu, N., Schaerf, A., Winter, F.
    \newblock Exact and metaheuristic approaches for unrelated parallel machine scheduling.
    \newblock {\em Journal of Scheduling}, 315-334, (2021).

    \bibitem{Oz18}	
    Ozturk, O., Chu, C.
    \newblock Exact and metaheuristic algorithms to minimize the total tardiness of cutting tool sharpening operations.
    \newblock {\em Expert Systems with Applications, 95}, 224-235, (2018).

    \bibitem{Par20}
    Parente, M., Figueira, G., Amorim, P., Marques, A., 
    \newblock Production scheduling in the context of Industry 4.0: review and trends.
    \newblock {\em International Journal of Production Research, 58:17}, 5401–5431, (2020).

    \bibitem{Pes22}
    Pessoa, A.A., Bulh\~oes, T., Nesello, V., Subramanian, A., 
    \newblock Exact Approaches for Single Machine Total Weighted Tardiness Batch Scheduling.
    \newblock {\em INFORMS Joural on Computing}, (2022).

    \bibitem{Peyro20}
    Peyro, L.F.,
    \newblock Models and an exact method for the Unrelated Parallel Machine scheduling problem with setups and resources.
    \newblock {\em Expert Systems with Applications: X, 5}, (2020).

    \bibitem{Pin20}
    Pinheiro, J.S.C.N., Aroyo, J.E.C., Fialho, L.B.,
    \newblock Scheduling unrelated parallel machines with family setups and resource constraints to minimize total tardiness.
    \newblock {\em GECCO '20: Proceedings of the 2020 Genetic and Evolutionary Computation Conference Companion}, 1409–1417, (2020).

    \bibitem{Pyomo}
    Pyomo - a Python-based Open-Source Software Package,
    \newblock {\em \url{https://www.pyomo.org/}}.

    \bibitem{Salah}
    Salah E. Elmaghraby, Park, S.H.,
    \newblock Scheduling Jobs on a Number of Identical Machines.
    \newblock {\em A I I E Transactions, 6:1}, 1-13 (1974).

    \bibitem{Yal00}	
    Yalaoui, F., Chu, C.
    \newblock Parallel Machine Scheduling to Minimize Total Completion Time with Release Dates Constraints.
    \newblock {\em IFAC Proceedings Volumes, 33:17}, 715-722, (2000).

    \bibitem{Yun22}	
    Yunusoglu, P., Yildiz, S.T.
    \newblock Constraint programming approach for multi-resource-constrained unrelated parallel machine scheduling problem with sequence-dependent setup times.
    \newblock {\em International Journal of Production Research, 60-7}, 2212–2229, (2022).
\end{thebibliography}
\end{document}